
\magnification=\magstep1
\input amstex
\input epsf.tex 
\documentstyle{amsppt}
\leftheadtext{P. Erd\H os, E. Makai, Jr., J. Pach}
\rightheadtext{Nearly equal distances in the plane, II}
\topmatter
\title Nearly equal distances in the plane, II
\vskip.3cm
\endtitle
\author P. Erd\H os, E. Makai, Jr.*, J. Pach{**} 
\vskip.3cm
{\centerline{\rm{Alfr\'ed R\'enyi Institute of Mathematics, HUN-REN}}}
{\centerline{\rm{H-1364 Budapest, P.O. Box 127, Hungary}}}
\vskip.1cm
{\centerline{\rm{http://www.renyi.hu/\~{}makai, http://www.renyi.hu/\~{}pach}}}
\vskip.1cm
{\centerline{\rm{E-mail: makai.endre\@renyi.hu,
pach.janos\@renyi.hu}}}
\vskip.1cm
{\centerline{\rm{ORCID ID: http://orcid.org/0000-0002-1423-8613,}}}
{\centerline{\rm{ORCID ID: http://orcid.org/0000-0002-2389-2035}}}
\endauthor

\thanks *Research (partially) supported by Hungarian National Foundation for 
Scientific Research, grant nos. 41, 1907, 2114, and ERC Advanced Grant
``GeoScape'' no. 882971.
\newline
{**}Research (partially) supported by Hungarian National Foundation for 
Scientific Research, grant no. 1907, and ERC Advanced Grant
``GeoScape'' no. 882971.
\endthanks


\abstract
Let $\{p_1, \ldots , p_n \} \subset {\Bbb{R}}^2$ be a separated point set, i.e.,
any two points have a distance at least $1$. Let $k \ge 1$ be an integer, and
$1 \le t_1 < \ldots < t_k$ be real numbers. Let $\delta > 0$. Suppose for all
$1 \le \ell (1) \le \ell (2) < \ell (3) \le k$ that
$|t_{\ell (3)} / (t_{\ell (1)} + t_{\ell (2)}) - 1| \ge
\delta
$. Then for $n \ge n_{k, \delta }$,
the number of pairs $\{ p_i,p_j\} $,
for which $d(p_i,p_j) \in [t_1, t_1 + 1] \cup \ldots \cup
[t_k, t_k + 1] 
$, is at most $n^2/4 + C_{k,\delta }n$. This is sharp, up
to the value of the constant $C_{k,\delta } > 0$.
\endabstract 
\endtopmatter
\document

%
%
%
%
%
%
%
%
%

\head
1. Introduction
\endhead

We call a finite set $\{ p_1, \ldots , p_n \} \subset {\Bbb{R}}^d$
{\it{separated}} if $\min _{i \ne j} d(p_i,p_j) \ge 1$, where $d(x,y)$ denotes
the distance of $x$ and $y$. We write $[x,y]$ for the closed
segment with
endpoints $x,y$, 
and ${\text{diam}}\,(\cdot )$
for the diameter of a set in ${\Bbb{R}}^d$.
We write $e_d := (0, \ldots , 0, 1) \in {\Bbb{R}}^d$.
In the paper $c$, $C$, $c_d$, etc. will denote
positive constants, which may have different values at
different occurrences.

It was a favourite question of the first author, what is the
maximal number of equal distances
determined by $n$ points in ${\Bbb{R}}^d$.
This question
induced much research in combinatorial geometry.
For the plane
he conjectured that
this maximum is $\Omega (n^{1 + c/\log \log n})$,
for some constant
$c > 0$, which is attained
for a ${\sqrt{n}} \times {\sqrt{n}}$ square grid (for $n$ a
square), cf. \cite{E1}.
The best known upper bound for this question in the plane
is $O(n^{4/3})$, cf. \cite{SST}.
For ${\Bbb{R}}^3$ the best known lower bound is $\Omega 
(n^{4/3} \log \log n)$, cf. \cite{E2}, and
the best known upper bound
is $O_{\varepsilon }(n^{3/2 - 1/394 + \varepsilon })$, cf.
\cite{Z}.
The case of ${\Bbb{R}}^4$ is completely solved in \cite{Br} and
\cite{vW}. Their results, taken together, say that
for $n \ge 5$
this maximum
is $\lfloor n^2/4 \rfloor + n$ if $8|n$ or $10|n$, and 
$\lfloor n^2/4 \rfloor + n - 1$ else.
Let $d \ge 6$ be even and let $n \ge n_d$. Then
this maximum is
$(n^2/2)(1 - 1/\lfloor d/2 \rfloor ) + n - O_d(1)$, and
for $2d|n$ it equals
$(n^2/2)(1 - 1/\lfloor d/2 \rfloor ) + n$,
cf. \cite{E3}. Again let $d \ge 6$ be even and let $n \ge n_d$.
Then the exact value of this maximum is determined by
\cite{S}.
For $d \ge 5$ odd and $n \ge n_d$ this maximum is $(n^2/2)
(1 - 1/\lfloor d/2 \rfloor ) + \Theta _d(n^{4/3})$, cf.
\cite{EP}.
For more information cf. \cite{PA}, Part II, Ch. 10.

About 1989 the first author observed that one can modify this question so that
the solution of the modified question is probably much simpler. Rather than
the number of equal distances, he asked for the number of almost equal
distances. Here several distances are {\it{almost equal}}
if they fall in a unit
interval. To avoid trivialities, when all distances are close to $0$, he
supposed that the system of $n$ points is separated.

In \cite{EMPS} it has been shown that for each $d \ge 2$
there exists a constant $n_d$, such that the following holds.
Let $n \ge n_d$, let
$t \ge 1$ and let
$\{ p_1, \ldots , p_n \} \subset {\Bbb{R}}^d$ be
a separated set. 
Then the number of
pairs
$\{ p_i,p_j \} $, for which $d(p_i,p_j) \in [t,t + 1]$, is at
most $(n^2/2)(1 - 1/d) - O_d(1)$. This bound
can be attained for any $n$.
We have equality for the set
$\{ s_{\mu} + \nu e_d \mid 1 \le \mu \le d, \,\, 1 \le \nu
\le n_{\mu } \} $.
Here
$s_1 \ldots s_d$ is a regular simplex in the $x_1 \ldots x_{d - 1}$-coordinate hyperplane, with vertices $s_{\mu }$, and
of edge length $t$. Moreover, for $1 \le \mu \le d$ we have
$\lfloor n/d \rfloor \le
n_{\mu } \le \lceil n/d \rceil $ and $\sum _{\mu = 1}^d
n_{\mu } = n$.
Further, we suppose that $t$ is sufficiently large
($t \ge c n^2$ suffices). 

In \cite{EMP1} the following has been proved. Let $k \ge 1$ be a
fixed integer, let $1 \le t_1 < \ldots < t_k$ be arbitrary, and let
$\{ p_1, \ldots , p_n \} \subset {\Bbb{R}}^2$ be any separated set. Then the
number of pairs $\{ p_i,p_j \} $, for which $d(p_i,p_j) \in
\cup _{\ell = 1}^k [t_{\ell }, t_{\ell } + 1] 
$, is at most $(n^2/2)\left( 1 - 1/(k + 1) \right) +
o_k(n^2)$. This estimate is sharp, up to the summand
$o_k(n^2)$. We have
equality (with $-O_k(1)$ rather than $o_k(n^2)$)
for the set
$\{ (\mu t, \nu) \mid 0 \le \mu \le k, \,\,1 \le \nu \le
n_{\mu } \} $.
Here $\lfloor n/(k + 1) \rfloor \le n_{\mu } \le \lceil
n/(k + 1) \rceil $ and
$\sum _{\mu = 0}^k n_{\mu } = n$, and 
$t_1 = t, \ldots , t_k = kt$, and $t$
is sufficiently large ($t \ge cn^2$ suffices).

Still we note that the above two results remain valid, if we
consider an interval
of the form $[t,t + c_d n^{1/d}]$,
or intervals of the form $[t_{\ell }, t_{\ell } +
c_{k, \varepsilon }{\sqrt{n}}]$. Here 
$c_d$ and
$c_{k, \varepsilon }$ are suitable constants.
In the second case $o_k(n^2)$ has to be substituted by
$\varepsilon n^2$, and we have to suppose that $n \ge n_{k,
\varepsilon }$.

A common generalization of these two theorems is the
following problem: generalize the result of \cite{EMP1} to
${\Bbb{R}}^d$.
This has been
treated for $k = 2$ in \cite{EMP2}. Their results have been
generalized for an
arbitrary $k$ by \cite{FK1} and \cite{FK2}. The last two
papers in a sense have settled the whole
problem, by reducing it to simpler problems (\cite{FK2}
is an extended abstract of \cite{FK1}). Also here one can
take, for $n \ge n_{d,k, \varepsilon }$,
intervals of the form $[t_{\ell }, t_{\ell } +
c_{d,k,\varepsilon } n^{1/d}]$.
Under these hypotheses, \cite{FK1}
determined for each $d$ and $k$ the maximum number of our point
pairs, up to a summand
$\varepsilon n^2$.
Moreover, for $d \ge d_k$ and
$n \ge n_{d,k}$ and for intervals $[t_{\ell }, t_{\ell } +
c_{d,k}n^{1/d}]$, \cite{FK1} determined
the exact maximum of the number of our point pairs. This
turned out to be a certain Tur\'an number.

\cite{PRV1} and \cite{PRV2} showed the foillowing. Let
$\{ p_1, \ldots , p_n \} \subset {\Bbb{R}}^d$, and let
$cn^2$ pairs $\{ p_i,
p_j \} $ have distances in some interval $[t,t+1]$. Then
for some $f_d(c) > 0$ we have
${\text{diam}}\, \{ p_1, \ldots , p_n \} \ge
f_d(c)n^{2/(d - 1)}$.
(Actually both papers showed $t \ge f_d(c)n^{2/(d - 1)}$.)
The order of magnitude is sharp.
An example for
$n = 2 m^{d - 1}$, where $m \ge 1$ is an integer, is
the set $\{ (i_1, \ldots , i_{d - 1}, x) \mid
1 \le i_1, \ldots, i_{d - 1} \le m {\text{ are integers,
and }} x \in \{ 0, c_d m^2 \} \} $.

Two, or $k$ nearly equal distances in arbitrary $n$-point
metric spaces have
been investigated in \cite{OP}. However, there nearness of
distances was
measured by the closeness of their quotients to $1$ (like
in \cite{EMP2}, Theorem 2). Moreover, their question was to
estimate from below the number $k$
of pairs of points $\{ p_i,q_i
\} $, with all quotients $d(p_{i(1)},q_{i(1)})/
d(p_{i(2)},q_{i(2)})$ in a small neighbourhood of $1$.

There arises the question, what can we tell for the planar case,
for a ``general'' system of values
$t_1, \ldots , t_k$. This is answered by our theorem, which
surprisingly asserts that the bound for any fixed $k$ is essentially the same
as for $k = 1$. Our result also can be interpreted as
``having many
distances almost equal to some $k$ distances
creates some order in the numbers $t_1, \ldots , t_k$''.


\head
2. The theorem and some problems
\endhead


\proclaim{Theorem}
Let $k \ge 2$ be a fixed integer, and
let $\delta \in (0,1)$ be fixed.
Then
there exist constants
$n_{k, \delta }, C_{k, \delta } > 0$, depending on $k$ and
$\delta $, with the following property. Suppose that
$n \ge n_{k, \delta }$, and 
$1 \le t_1 < \ldots < t_k$ satisfy the following hypothesis:
$$
\cases
{\text{for any }} 1 \le \ell (1) \le \ell (2) < \ell (3) \le k
{\text{ we
have }}
\\
t_{\ell (3)} \not\in
[(1 - \delta) (t_{\ell (1)} + t_{\ell (2)}) ,
t_{\ell (1)} + t_{\ell (2)} + 2] .
\endcases
$$
Then for any separated set $\{ p_1, \ldots , p_n \} \subset {\Bbb{R}}^2$, the
number of pairs $\{ p_i,p_j \} $ for which $d(p_i,p_j) \in
\cup _{\ell = 1}^k
[t_{\ell }, t_{\ell } + 1]
$, is at most $n^2/4 + C_{k, \delta }n$.

Moreover, the same inequality holds for intervals
$[t_{\ell }, t_{\ell }
+ \alpha ]$ rather than $[t_{\ell }, t_{\ell } + 1]$, for 
any fixed constant $\alpha $. However, then we have
$n_{k, \delta , \alpha }$ and $C_{k, \delta , \alpha }$,
rather than $n_{k, \delta }$ and $C_{k, \delta }$.

The estimate of the theorem is essentially sharp. Even,
for any fixed $\varepsilon > 0$, for each $k$ and each
sufficiently small $\delta > 0$ and each $n$
we have the following.
There exist examples
with intervals $[t_{\ell }, t_{\ell } + \varepsilon ]$
rather than
$[t_{\ell }, t_{\ell } + 1]$,
for which the number of the above pairs is
$n^2/4 + (k - 1)n - O_k(1)$.
Moreover,
this holds for all sufficiently large $t_k$, depending on
$k$, $\varepsilon $ and $n$ ($t_k \ge \max \{ 3^{k - 1},
cn^2/\varepsilon \} $ suffices).
\endproclaim


{\bf{Remark 1.}}
The simplified statement of the Theorem in the abstract
follows from this form of the Theorem, provided $t_{\ell (1)}
$ is sufficiently large. If however, $t_{\ell (1)} \ge t_1$ is
bounded, then {\bf{1}} of the proof of the Theorem shows the
simplified statement of the Theorem in the abstract.

\vskip.1cm


{\bf{Remark 2.}}
If we allow $t_{\ell (3)} = t_{\ell (1)} + t_{\ell (2)}$
for $1 \le \ell (1) \le \ell (2) < \ell (3) \le k$, then
we may have $\lfloor n^2/3 \rfloor $
pairs with distances in $\cup _{\ell = 1}^k [t_{\ell },
t_{\ell } + 1] 
$. An example is $\{ (0,1), \ldots ,
(0,n_1),$
$(t_{\ell (1)}, 1), \ldots , (t_{\ell (1)},n_2),
(t_{\ell(3)},1),
\ldots (t_{\ell (3)},n_3) \}
$. Here $\lfloor n/3 \rfloor \le n_{\mu } 
\le \lceil n/3 \rceil $ and $\sum _{\mu = 1}^3 n_{\mu }
= n$, and $t_{\ell (1)}$ is sufficiently large
($t_{\ell (1)} \ge cn^2$
suffices). Actually an analogous example
exists if
$|t_{\ell (3)} - t_{\ell (1)} - t_{\ell (2)}|$ is at most
some fixed number in $(0,1)$.


\vskip.1cm

{\bf{Problem 1.}}
Can one prove a similar theorem under the following
weaker hypothesis: for 
$1 \le \ell (1) \le \ell (2) < \ell (3) \le k$
we have $t_{\ell (3)}
\not\in [t_{\ell (1)} + t_{\ell (2)} -
C (t_{\ell (1)} + t_{\ell (2)})^c,
t_{\ell (1)} + t_{\ell (2)} + 2]$, for some $c \in [0,1)$ and
$C > 0$?


\vskip.1cm

{\bf{Problem 2.}}
Analogously like in 
\cite{EMPS}, \cite{EMP1}, \cite{EMP2}, \cite{FK1}, is it
possible to extend our Theorem to intervals of the form
$[t, t + f(n)]$, for some $f(n)$ tending to $\infty $
for $n \to \infty $?


\vskip.1cm

{\bf{Problem 3.}}
For $k = 2$ our theorem settles ``everything''. For $t_2/t_1 = 2$ we have the
essentially sharp bound $n^2/3 + \varepsilon n^2$ (\cite{EMP1}). However, for
$t_2/t_1$ far from $2$ we have a bound $n^2/4 + Cn$. Now let $k \ge 3$.
Suppose that we want to exclude only the systems $(t_1, \ldots , t_k) = (t,
\ldots , kt)$, and those ``close'' to them. (E.g., those which satisfy
$|t_{\ell }/t_1 - \ell | \le \delta $ for all $2 \le \ell
\le k$.)
Then conjecturably, for $n$ sufficiently large, $n \ge n_k$,
say, the maximum
number of pairs of points with distances in
$\cup _{\ell = 1}^k [t_{\ell },t_{\ell } + 1]$
is obtained as follows. We let $\{ p_1, \ldots , p_n \} :=
\{ (\mu t, \nu) \mid 1 \le \mu \le k,\,\,1 \le \nu \le
n_{\mu } \} $. Here for $1 \le \mu \le k$ we have
$ 
\lfloor n/k \rfloor \le n_{\mu } \le \lceil n/k \rceil $ and
$\sum _{\mu = 1}^k n_{\mu } = n$. Moreover, we take
$t_1 = 1$, $t_2 = t, \ldots $, $t_k = (k-1)t$.
The number of pairs in
question for this example is $(n^2/2)(1 - 1/k) + 2n - O_k(1)$.
Actually this was the original question of the first author. Can one formulate
analogous questions, with suitable hypotheses,
for which the maximum number of pairs with
distances in $\cup _{\ell = 1}^k [t_{\ell },t_{\ell } + 1]$
is conjecturably
$(n^2/2)(1 - 1/s) + o_k(n^2)$, for
$3 \le s \le k - 1$? (Heuristically: if the number of our
distances almost equal to some $k$ numbers $t_1, \ldots , t_k$
is $(n^2/2)(1 - 1/s) + o_k(n^2)$, then with increasing $s$ we
have ever more structure in $\{ t_1, \ldots , t_k \} $.)


\vskip.1cm

{\bf{Problem 4.}} Can one generalize our theorem to 
${\Bbb{R}}^d$, in the following form. Suppose that
there is no ``close to degenerate''
simplex with all edge lengths some $t_{\ell }$'s
(degenerate meaning
that it lies in a hyperplane). 
Then the number of pairs $\{ p_i, p_j \} $, for which
$d(p_i,p_j)
\in \cup _{\ell = 1}^k [t_{\ell },t_{\ell } + 1]$, is at most
$(n^2/2)(1 - 1/d) + C_{d,k} n$. Here of course $C_{d,k}$
should depend also
on the ``measure of non-degenerateness'' (this is $\delta $
in our Theorem). A possible measure of degeneracy could be
$(d, \alpha )$-Flatness of point-sets from \cite{FK1}.

This estimate, if true, would be essentially sharp.
In fact, the example from
\cite{EMPS}, cited in the introduction, with number of pairs
in question being
$(n^2/2)(1 - 1/d) - O_d(1)$,
works for any $k \ge 1$.
(Moreover, we could choose $t_1 := 1$, $t_2 := 3$, $\ldots $,
$t_{k - 1} := 3^{k - 2}$, and $t_k := t \ge \max
\{ 3^{k - 1},cn^2 \} $, and any $\delta \le 1/3$,
provided that with edge lengths
these $t_{\ell }$'s there is no ``close to degenerate''
simplex. Then we would have even
$(n^2/2)(1 - 1/d) + 2(k - 1)n - O_{d,k}(1)$ pairs in
question.)

Now suppose that
there is a degenerate simplex $s_1 \ldots s_{d + 1}$
with vertices $s_{\mu }$ and
with all edge lengths some $t_{\ell }$'s, all these
$t_{\ell }$'s being
at least $cn^2$. Let $s_1 \ldots s_{d + 1}$ lie
in the $x_1 \ldots x_{d - 1}$-coordinate hyperplane.
Then there
exist examples for which the number of the pairs in question is
$(n^2/2)\left( 1 - 1/(d + 1) \right) 
- O_{d,k}(1)$.
Our example is the $d$-dimensional
generalization of the example in Remark 2.
Suppose that for $1 \le \mu \le d + 1$ we have
$\lfloor n/(d + 1)
\rfloor \le n_{\mu } \le \lceil n/(d + 1) \rceil $ and
$\sum _{\mu = 1}^{d + 1} n_{\mu } = n$.
Then our example is
$\{ s_{\mu } + \nu e_d \mid 1 \le \mu \le d + 1,\,\,1 \le
\nu \le n_{\mu } \} $.

Cf. also Remark 3 in the end of \S 3.


\head
3. The proof of the theorem
\endhead


\demo{Proof}
{\bf{1.}}
First we suppose that
$$
t_1 \le C, {\text{ where }} C {\text{ is some (large)
constant.}}
\tag 1
$$

Then, like in \cite{EMP1}, we use induction.
The base of induction
is $k = 1$, for which the statement holds by \cite{EMPS},
cited in the introduction.
For $k \ge 2$ the hypothesis of the theorem is satisfied for
$1 \le t_2 < \ldots <t_k$. Hence
the number of pairs
$\{ p_i,p_j \} $, such that $d(p_i,p_j) \in \cup _{ \ell = 2}^k
[t_{\ell }, t_{\ell } + 1]
$, is at most $n^2/4 + C_{k - 1, \delta } n $. The number
of pairs $\{ p_i,p_j \} $ such that $d(p_i,p_j) \in
[t_1, t_1 + 1]$, and $p_i$
is fixed, is at most the quotient of the
area of a circular ring of radii $t_1 - 1/2$ and
$t_1 + 3/2$, and of $\pi /4$.
Namely, the open circles of radius $1/2$ centred at these
$p_j$'s are disjoint, and lie in this circular ring.
The above quotient is $2 \pi (2t_1 + 1) / (\pi /4) \le 8(2C + 1)$.
Thus the number of all pairs 
$\{ p_i,p_j \} $, such that $d(p_i,p_j) \in  [t_1, t_1 + 1]$,
is at most
$8(2C + 1)n$. Hence we can choose $C_{k, \delta } := C_{k - 1, \delta } +
8(2C + 1)$.

{\bf{2.}}
Therefore we may suppose that
$$
t_1 \ge C, {\text{ where }} C {\text{ is a suitable (large)
constant.}}
\tag 2
$$

We consider the graph $G$
with vertices $p_1, \ldots , p_n$, and edges those pairs
$\{ p_i, p_j \} $, for which $d(p_i, p_j) \in
\cup _{\ell = 1}^k [t_{\ell },t_{\ell } + 1]
$. Suppose, in contradiction to the statement of the theorem, that
the number of edges of $G$ is greater than $n^2/4 + C_{k, \delta }n$.
Then, for any preassigned number $N$,
for $C_{k, \delta }$ sufficiently large, 
depending on $N$, we have the following. Our
graph
$G$ contains a subgraph $K(1, N, N)$, cf. \cite{Bo1}, Cor. 4.7,
or \cite{Bo2}, Th. 1.5.2.
Here
$K(a,b,d)$ denotes the complete $3$-partite graph with colour classes $A,B$ and $D$,
of sizes $a,b$ and $d$, resp. Now $A,B,D \subset
\{ p_1, \ldots , p_n \} $, and $|A| = 1$
and
$|B| = |D| = N$. Let $\{ x \} = A$ and $y \in B$ and $z \in D$. Then the
triangle $xyz$ has all side-lengths in $\cup _{\ell = 1}^k
[t_{\ell },t_{\ell } + 1]$.

We let correspond to each side of the
triangle $xyz$ the
smallest $\ell $ such that 
the length of this side lies in $[t_{\ell },t_{\ell } + 1]$.
This $\ell $
will be denoted, for the sides $xy$, $yz$ and $zx$, by
$\ell (xy)$, $\ell (yz)$ and $\ell (zx)$, resp.

For $y \in B$ the number
$\ell (xy)$ can assume $k$ values. Hence,
for some $B_1 \subset B$ with $|B_1| = \lceil N/k \rceil $,
for $y \in B_1$ the number
$\ell (xy)$ is independent of
$y$. Similarly, 
for some $D_1 \subset D$ with $|D_1| = \lceil N/k \rceil $,
for $z \in D_1$ the number
$\ell (zx)$ is independent of
$z$. Observe that in our theorem $k$ is fixed. 

For $y \in B_1$ and $z \in D_1$ the number $\ell (yz)$ 
can assume $k$ values.
By a variant of Ramsey's
theorem, cf. \cite{PA}, Cor. 9.15 (stated there for two
colours only, but its proof given there
works for any fixed number of
colours -- alternatively, the statement for two colours
implies the statement for any fixed number of colours),
we have the following.
There exist
$B_2 \subset B_1$ and $D_2 \subset D_1$, such that $|B_2| =
|D_2| = M$, where 
also $M$ is arbitrarily large provided $N$ is arbitrarily
large,
such that the following holds.
For all pairs $\{ y, z \} $, where $y \in B_2$ and
$z \in D_2$, the number $\ell (yz)$ is independent of $y,z$.

Hence
$$
\cases
{\text{for all triangles }}xyz, {\text{ where }} \{ x \} = A,
\,\,y \in B_2 {\text{ and }} z \in D_2,
\\
{\text{we have that
}} \ell (xy),\,\, \ell (yz) {\text{ and }} \ell (zx)
{\text{ are
independent of }} y,z.
\endcases
\tag 3
$$

Let
$$
\{ x,y,z \} = \{ u,v,w \} , {\text{ where }}
\ell (1) := \ell (uv) \le \ell (2) := \ell (vw) \le \ell (3)
:= \ell (wu) .
\tag 4
$$
There are two cases.
$$
{\text{Either (I): }} \ell (1) \le \ell (2) < \ell (3),
\tag 5
$$
$$
{\text{or (II): }} \ell (1) \le \ell (2) = \ell (3).
\tag 6
$$

{\bf{3.}}
In case (I) we have
$$
t_{\ell (3)} \le d(w,u) \le d(u,v) + d(v,w) \le
t_{\ell (1)} + t_{\ell (2)} + 2 .
\tag 7
$$ 
This implies by the hypothesis of the theorem that
$$
t_{\ell (3)} \le (1 - \delta )
\left( t_{i(1)} + t_{i(2)} \right) .
\tag 8
$$
Since by \thetag{2}
$$
C \le t_1 \le t_{\ell (1)} \le t_{\ell (2)} < t_{\ell (3)},
\tag 9
$$
by choosing $C$ sufficiently large
($C \ge c/ \delta $ suffices), we will have
also
$$
d(w,u) \le t_{\ell (3)} + 1
\le (1 - \delta / 2) \left( t_{\ell (1)} + t_{\ell (2)}
\right) 
\le  (1 - \delta / 2) \left( d(u,v) + d(v,w) \right) .
\tag 10
$$
Moreover,
by the definition of $\ell (uv)$ etc., we have
$$
\max \{ d(u,v), d(v,w) \} < d(w,u) .
\tag 11
$$
By \thetag{10}, the
triangle $uvw$ cannot degenerate to the doubly counted segment
$[w,u]$, therefore it is not degenerate.

We claim
$$
\angle uvw > \max \{ \angle vwu, \angle wuv \} \ge \min
\{ \angle vwu, \angle wuv \} \ge \delta _1 ,
\tag 12
$$
where $\delta _1 = \delta _1 (\delta ) > 0$. The first
inequality follows from \thetag{11}. For the last inequality
we  prove
$\angle vwu \ge \delta _1$; the proof of $\angle wuv \ge \delta _1$ is the
same. 

In
fact, let the sides of $\angle vwu$ be fixed, and let
$d(v,w) \in [0, d(w,u)]$. Then, by a triangle inequality,
the maximum of $d(u,v) + d(v,w)$ is attained for
$d(v,w) = d(w,u)$. The value of this maximum is
$d(w,u) + 2d(w,u) \sin (\angle vwu /2)$. Hence, by
\thetag{10}, 
$$
\cases
d(w,u) \le (1 - \delta /2) \left( d(u,v) + d(v,w) \right) \le
\\
(1 - \delta /2) \left( d(w,u) + 2d(w,u) \sin (\angle vwu /2)
\right) .
\endcases
\tag 13
$$
Comparing the first and last expressions in \thetag{13}, we get
$$
\angle vwu \ge
2 \arcsin \left( \delta /(4 - 2 \delta ) \right) =:
\delta _1 \sim \delta /2 ,
{\text{and similarly }} \angle wuv \ge \delta _1 ,
\tag 14
$$
which ends the proof of \thetag{12}. Then 
\thetag{12} implies that
$$
\max \{ \angle vwu , \angle wuv \} < \angle uvw = \pi -
\angle vwu  - \angle wuv \le
\pi - 2 \delta _1 =: \pi - \delta _2 ,
\tag 15
$$
where $\delta _2  = \delta _2 (\delta )$.

The above considerations apply to all of
the triangles
$xyz$ from \thetag{3}. Let us choose one of them, $xy_0z_0$,
say.
Then both $y_0$ and $z_0$ can be replaced by any of $M$
points $y_i \in B_2$ and $z_i \in D_2$, say,
so that we obtain a triangle $xy_iz_i$ satisfying \thetag{3}.

Now we will
only use the fact that $y_0$ can be replaced by any of $M$
points $y_i \in B_2$,
so that we obtain a triangle $xy_iz_0$ satisfying \thetag{3}.
Then by \thetag{4}
all distances $d(x,y_i)$ belong to one of the intervals
$[t_{\ell (1)},
t_{\ell (1)} + 1]$,  $[t_{\ell (2)},t_{\ell (2)} + 1]$ and
$[t_{\ell (3)},t_{\ell (3)} + 1]$. Similarly,
all distances $d(y_i,z_0)$ belong to one of these
three intervals (possibly different from the above one).

Above we have chosen some $y_i \in B_2$, namely, $y_0$. Then
for each
$y_i \in B_2$ we have 
$|d(x,y_i) -
d(x,y_0)| \le 1$, i.e., $d(x,y_i) \in [d(x,y_0) - 1, d(x,y_0) + 1]$.
Similarly $d(y_i,z_0) \in [d(y_0,z_0) - 1, d(y_0,z_0)
+ 1]$. Moreover, by \thetag{12} and \thetag{15} we
have $\angle xy_0z_0 \in [\delta _1, \pi - \delta _2]$.

That is, all of these $M$ points $y_i$ lie in the
intersection of two circular rings, of centres $x$ and
$z_0$. These circular rings have inner radii $d(x,y_0)
- 1 \ge C - 1$
and $d(y_0,z_0) - 1 \ge C - 1$,
and outer radii $d(x,y_0) + 1$ and $d(y_0,z_0) + 1$, resp.
Since even $C - 1$ is large, this intersection is the union
of two
``curvilinear almost rhombs'', symmetric to the line $z_0x$.
Moreover, $\angle xy_0z_0 \in [\delta _1, \pi - \delta _2]$.

Any of these two ``curvilinear almost rhombs''
can be included in an intersection of two parallel strips
of width $3$, which is a rhomb $R$. These strips are orthogonal to
the vectors $y_0 - x$ and $y_0 - z_0$, resp.
The boundary lines of the first strip intersect the halfline
from $x$, passing through $y_0$, at points with distances
$d(x,y_0) + 1$ and $d(x,y_0) - 2$ from $x$. Similarly, the
boundary lines of the second strip intersect the halfline
from $z_0$, passing through $y_0$, at points with distances
$d(y_0,z_0) + 1$ and $d(y_0,z_0) - 2$ from $z_0$.
This statement is very intuitive; a formal proof can be found in
the first arxiv variant of \cite{EMP2}, pp. 9-13. Then all $M$
points $y_i \in B_2$ lie in the union of $R$ and its axially
symmetric image w.r.t. the line $z_0x$.

Consider the rhomb $R_1$, whose side lines are
obtained from those of the rhomb $R$ by translating them
outwards (from the centre of $R$) through a distance $1/2$.
The height of $R_1$ is $4$.
By $\angle xy_0z_0 \in [\delta _1, \pi - \delta _2]$ the
angles of 
$R_1$ are bounded away from $0$ and $\pi $. Hence its area, which is $4^2$
divided by the sine of one of its angles, i.e., by $\sin \angle xy_0z_0$,
is bounded above by
some function of $\delta $. Observe that in our Theorem
$\delta $ is fixed.

All $M$ disjoint open circles of radius $1/2$, of
centres $y_i \in B_2$,
are
contained in the union of $R_1$ and its axially symmetric
image w.r.t.
the line
$z_0x$. Hence
by area considerations $M$ must be bounded. However,
$M$ can be arbitrarily large, a contradiction.

{\bf{4.}}
In case (II) (cf. \thetag{6})
we choose one of the triangles, $xy_0z_0$, say.
By
\thetag{4} and \thetag{6} we have
$\{ x,y_0,z_0 \} = \{ u_0,v_0,w_0 \} $, where
$$
\ell (u_0,v_0) = \ell (1) \le \ell (v_0,w_0) = \ell (2) =
\ell (w_0,u_0) = \ell (3).
\tag 16
$$
Then one of $u_0$ and $v_0$ is
different from $x$, hence belongs to $B_2$ or to $D_2$.
Let, e.g., $u_0 \in B_2$. Then
$u_0$ can be replaced by any of $M$ points $y_i \in B_2$,
so that we obtain a triangle $xy_iz_0$ satisfying \thetag{3}.

Now we claim that
$\angle w_0u_0v_0$ lies in a small neighbourhood of the angular interval
$[\pi /3, \pi /2]$. In fact, let us write $a := d(w_0,u_0)$,
$b := d(u_0,v_0)$ and $c := d(v_0,w_0)$. Then
we have for the angle $\gamma $ of the triangle $u_0v_0w_0$,
opposite to its side $c$, that
$$
\cases
\cos \gamma \le [(t_{\ell (2)} + 1)^2 + (t_{\ell (1)} + 1)^2
-
t_{\ell (2)}^2]/[2t_{\ell (2)} t_{\ell (1)}] =
\\
t_{\ell (1)}/(2t_{\ell (2)}) + 1/t_{\ell (1)}
+ 1/t_{\ell (2)} +
1/(t_{\ell (1)}t_{\ell (2)}) \le 1/2 + 2/C + 1/C^2,
\endcases
\tag 17
$$
where $C$ can be anyhow large. Therefore $\gamma $ is at least a value
close to $\pi /3$.

On the  other hand,
$$
\cos \gamma \ge [t_{\ell (2)}^2 + t_{\ell (1)}^2 -
(t_{\ell (2)} + 1)^2]/(2ab).
\tag 18
$$
If the numerator of \thetag{18} is nonnegative, then
$\gamma \le \pi /2$.
If this numerator is negative, then the right hand side of
\thetag{18}
can be estimated further from below by
$$
\cases
[t_{\ell (1)}^2 - 2t_{\ell (2)} -
1]/(2t_{\ell (1)}t_{\ell (2)}) =
t_{\ell (1)}/(2t_{\ell (2)})
\\
- 1/ t_{\ell (1)} -
1/(2t_{\ell (1)}t_{\ell (2)})
\ge 0 - 1/C - 1/(2C^2),
\endcases
\tag 19
$$
and then $\gamma $ is at most a value close to $\pi /2$.

Summing up: $\gamma $ is bounded away both from $0$ and from $\pi $.

Then using $|B_2| = M$, which can be arbitrarily large, and the area
considerations from the end of
{\bf{3}}, we get a contradiction.

{\bf{5.}}
Now we prove the assertion with $\alpha $, where we may
suppose $\alpha > 1$. In {\bf{1}} the base of induction is
the statement from \cite{EMPS}, where we can have intervals
of the form $[t, t + c {\sqrt{n}}]$. So we may take
intervals of the form $[t, t + \alpha ]$, for $n \ge
n_{k, \delta, \alpha }$. After this point,
let us diminish $ \{ p_1, \ldots , p_n
\} $ in ratio $1/\alpha $. Then rather than $[t_{\ell },
t_{\ell } + \alpha ]$
we will have intervals
$[t_{\ell }/ \alpha , t_{\ell }/ \alpha + 1]$, and rather than $t_1 \ge 1$
we will have $t_1/ \alpha \ge
1/\alpha $. Furthermore, all the steps of the
above proof apply for this diminished point system, except
that the minimal distance is at least $1/ \alpha $. Then
we will
have disjoint open circles of radius $1/(2 \alpha )$, where
$\alpha $ is a constant. The area considerations work for
these disjoint open circles as well, in {\bf{1}}, in the end
of {\bf{3}} and in the end of {\bf{4}}. 

{\bf{6.}}
The asserted essential sharpness of the estimate is shown by
the example which is the planar case of the example
from \cite{EMPS}, cited in the introduction. Here we
choose $t_1,
\ldots , t_k$ as in the second paragraph of Problem 4.
The statement about essential sharpness
with the intervals $[t_{\ell }, t_{\ell } + \varepsilon ]$
follows from this by Pythagoras' theorem. 
$\blacksquare $
\enddemo


{\bf{Remark 3.}}
Concerning Problem 4,
unfortunately the proof of our Theorem
does not carry over even to $d = 3$ and $k = 2$. E.g., we can
have tetrahedra with two different edge lengths $t_1 < t_2$,
two opposite
edges of length $t_1$, the remaining ones of length $t_2$,
where $t_2/t_1$ is large. For these for each vertex the
absolute value of the
determinant formed by the unit vectors pointing from that
vertex to all the other vertices is small. (However, all other
tetrahedra with these edge lengths can probably
be treated as in the proof of our Theorem.)
However,
using $(d, \alpha )$-Flatness of point-sets
from \cite{FK1}, this
example is excluded. So this gives a chance to apply
$(d, \alpha )$-Flatness of point-sets
as a measure of non-degeneracy in
${\Bbb{R}}^d$.


\end 


(E.g., the quotient of
the volumes of the inball and of the simplex is bounded
below. In the plane this means that the minimal
angle is bounded below.)


and for $n \ge n_{d,k}$ (where $n_{d,k}$
is some constant)


and for $n \ge n_{d,k}$


for each $k$ there exists a constant
$n_k$, such that the following holds.


and all $n \ge n_k$,


{\bf{Remark 2.}} Concerning Problem 4, 
the combinatorial part of the proof of our Theorem goes
through also for ${\Bbb{R}}^d$.
This is clear for the Simonovits theorem. For the
Ramsey-type theorem suppose that
$d$ and $k$ are fixed,
and let $K_{1,N, \ldots , N}$ be a complete $d$-partite graph,
with $N$ sufficiently large. Then for
any $k$-colouring of the edges of the graph
$K_{1,N, \ldots , N}$, this graph
contains a subgraph
$K_{1,M, \ldots , M}$ with $M$ large, 
with the $d$ colour classes
of $K_{1,M, \ldots , M}$ contained in the respective $d$ colour
classes of $K_{1,N, \ldots , N}$. Moreover, $K_{1,M, \ldots , M}$
satisfies that the colour of an
edge of it only depends on the colour
classes of $K_{1,N, \ldots , N}$ which the edge endpoints
belong to.
This can be
proved by applying this statement for the bipartite case
(used in the proof of our Theorem),
successively for all $d \choose 2$
pairs of the colour classes of $K_{N, \ldots , N}$.

Moreover, we may assume $t_1 \ge
C$, where $C$ is a suitable (large) constant.

There arises the
question, how to measure ``non-degenerateness'' of a
simplex. Surely we have to exclude all degenerate simplices
$S$
with all edge-lengths some $t_i$'s (cf. Problem 4),
but also those ``close'' to them.  
One ought to
look for such a measure of ``non-degenerateness'' of a
simplex, that if this measure is at least some $\delta > 0$,
then at least for two vertices of the simplex we should have
that the absolute value of the determinant with columns
the unit vectors pointing from this vertex to all other
vertices should be at least some $f(\delta ) > 0$.
Then one
could repeat the arguments about boundedness of the volume
of the higher dimensional analogs of the ``curvilinear almost
rhombs'' (intersections of $d$ spherical shells of width
$2$ and of large inner radii) and finish the proof.


let $s_1 \ldots s_d$ be a
regular simplex
in the $x_1 \ldots x_{d - 1}$-coordinate
hyperplane, with vertices
$s_{\ell }$, and edge length $t \ge cn^2 
$. Suppose that for $1 \le i \le d$ we have $\lfloor n/d \rfloor
\le n_i \le \lceil n/d \rceil $ and $\sum _{i = 1}^d n_i
= n$.
Then an example is $\{ s_1 + e_d, \ldots , s_1 + n_1 e_d,
\ldots , s_d + e_d, \ldots ,s_d + n_d e_d \} $. 


$\{ s_1 + e_d, \ldots , s_1 + n_1 e_d, \ldots , s_{d + 1} +
e_d, \ldots ,s_{d + 1} + n_{d + 1} e_d \} $.


Moreover,
their boundary
lines intersect the halflines from $x$ and $z_0$, resp.,
passing through
$y_0$, at points with distances $d(x,y_0) + 1$ and $d(x,y_0)
- 2$
from $x$, and at points with distances
$d(y_0,z_0) + 1$ and $d(y_0,z_0) - 2$ from $z_0$,
resp.
